\documentclass[11pt,twoside,center]{article}

\usepackage{amsmath,amssymb}

\newtheorem{definition}{Definition}
\newtheorem{theorem}{Theorem}
\newtheorem{example}{Example}

\newtheorem{remark}{Remark}
\newtheorem{proposition}{Proposition}

\newenvironment{proof}[1][Proof]{\noindent\textbf{#1.} }

\title{\textsc{Dichotomies for evolution equations in Banach spaces }}

\author{Codru\c{t}a Stoica}
\date{~}
\begin{document}

\maketitle

\pagestyle{myheadings} \markboth{Codru\c{t}a Stoica}{Dichotomies
for evolution equations}

\begin{abstract}
The aim of this paper is to emphasize various concepts of
dichotomies for evolution equations in Banach spaces, due to the
important role they play in the approach of stable, instable and
central manifolds. The asymptotic properties of the solutions of
the evolution equations are studied by means of the asymptotic
behaviors for skew-evolution semiflows.


\end{abstract}

\section{Preliminaries}

Recently, the important progress made in the study of evolution
equations had a master role in the developing of a vast
literature, concerning mostly the asymptotic properties of linear
operators semigroups, evolution operators or skew-product
semiflows.

In this paper, the study is led throughout the notion of
skew-evolution semiflow on Banach spaces, defined by means of an
evolution semiflow and an evolution cocycle. As the
skew-evolutions semiflows reveal themselves to be generalizations
of evolution operators and skew-product semiflows, they are
appropriate to study the asymptotic properties of the solutions of
evolution equations having the form
\[
\left\{
\begin{array}{l}
\dot{u}(t)=A(t)u(t), \ t>t_{0}\geq 0 \\
u(t_{0})=u_{0},%
\end{array}%
\right.
\]
where $A:\mathbf{R}\rightarrow \mathcal{B}(V)$ is an operator,
$\textrm{Dom}A(t)\subset V$, $u_{0}\in \textrm{Dom}A(t_{0})$.

The fact that a skew-evolution semiflow depends on three variables
$t$, $t_{0}$ and $x$, while the classic concept of cocycle depends
only on $t$ and $x$, justifies the study of asymptotic behaviors
in a nonuniform setting (relative to the third variable $t_{0}$)
for skew-evolution semiflows.

The basic concepts of asymptotic properties, such as stability,
instability and dichotomy, that appear in the theory of dynamical
systems, play an important role in the study of stable, instable
and central manifolds. We intend to define and exemplify various
concepts of dichotomies, as uniform exponential dichotomy,
Barreira-Valls exponential dichotomy, exponential dichotomy,
uniform polynomial dichotomy, Barreira-Valls polynomial dichotomy,
polynomial dichotomy and to emphasize connections between them. We
have thus considered generalizations of some asymptotic properties
for evolution equations, defined by L. Barreira and C. Valls in
\cite{BaVa_LNM}. Characterizations for the asymptotic properties
in a nonuniform setting are also proved.

Some of the original results concerning the properties of
stability and instability for skew-evolution semiflows      were
published in \cite{StMe_NA} and \cite{MeSt_ICNODEA07}.

The exponential dichotomy for evolution equations is one of the
mathematical domains with an impressive development due to its
role in describing several types of differential equations. Its
study led to an extended literature, which begins with the
interesting results due to O. Perron in \cite{Pe_MZ}. The ideas
were continued by J.L. Massera and J.J. Sch\"{a}ffer in
\cite{MaSc_PAM}, with extensions in the infinite dimensional case
accomplished by J.L. Dalecki\u{i} and M.G. Kre\u{i}n in
\cite{DaKr_AMS} and A. Pazy in \cite{Pa_SV}, respectively R.J.
Sacker and G.R. Sell in \cite{sacsel}. Diverse and important
concepts of dichotomy were introduced and studied by S.N. Chow and
H. Leiva in \cite{ChLe_JDE} and by W.A. Coppel in \cite{Co_LNM}.

Some asymptotic behaviors for evolution families were given in the
nonuniform case in \cite{MeSaSa_MIR} by M. Megan, A.L. Sasu and B.
Sasu. The study of the nonuniform exponential dichotomy for
evolution families was considered by P. Preda and M. Megan in
\cite{PrMe_BAMS}.

The property of exponential dichotomy for the case of
skew-evolution semiflows is treated in \cite{MeSt_TP09} and
\cite{StMe_OT}.

\section{Notations. Definitions. Examples}

Let us consider a metric space $(X,d)$, a Banach space $V$,
$V^{*}$ its topological dual and $\mathcal{B}(V)$ the space of all
bounded linear operators from $V$ into itself. $I$ is the identity
operator on $V$. We denote $T =\left\{(t,t_{0})\in \mathbf{R}^{2},
\ t\geq t_{0}\geq 0\right\}$ and $Y=X\times V$.

\begin{definition}\rm\label{def_sfl_ev}
A mapping $\varphi: T\times  X\rightarrow  X$ is called
\emph{evolution semiflow} on $ X$ if following relations hold:

$(s_{1})$ $\varphi(t,t,x)=x, \ \forall (t,x)\in
\mathbf{R}_{+}\times X$;

$(s_{2})$ $\varphi(t,s,\varphi(s,t_{0},x))=\varphi(t,t_{0},x),
\forall (t,s),(s,t_{0})\in T, x\in X$.
\end{definition}

\begin{definition}\rm\label{def_aplcoc_ev}
A mapping $\Phi: T\times  X\rightarrow \mathcal{B}(V)$ is called
\emph{evolution cocycle} over an evolution semiflow $\varphi$ if:

$(c_{1})$ $\Phi(t,t,x)=I$, $\forall (t,x)\in \mathbf{R}_{+}\times
X$;

$(c_{2})$
$\Phi(t,s,\varphi(s,t_{0},x))\Phi(s,t_{0},x)=\Phi(t,t_{0},x),\forall
(t,s),(s,t_{0})\in T, x\in X$.
\end{definition}

\begin{definition}\rm\label{def_coc_ev_1}
The mapping $C: T\times Y\rightarrow Y$ defined by the relation
$$C(t,s,x,v)=(\varphi(t,s,x),\Phi(t,s,x)v),$$ where $\Phi$ is an
evolution cocycle over an evolution semiflow $\varphi$, is called
\emph{skew-evolution semiflow} on $Y$.
\end{definition}

\begin{example}\rm\label{ex_ses}
We denote by
$\mathcal{C}=\mathcal{C}(\mathbf{R}_{+},\mathbf{R}_{+})$ the set
of all continuous functions $x:\mathbf{R}_{+}\rightarrow
\mathbf{R}_{+}$, endowed with the topology of uniform convergence
on compact subsets of $\mathbf{R}_{+}$, metrizable by means of the
distance
\[
d(x,y)=\sum_{n=1}^{\infty}\frac{1}{2^{n}}\frac{d_{n}(x,y)}{1+d_{n}(x,y)},
\ \textrm{where} \ d_{n}(x,y)=
\sup\limits_{t\in[0,n]}{|x(t)-y(t)|}.
\]
If $x\in \mathcal{C}$, then, for all $t\in \mathbf{R}_{+}$, we
denote $x_{t}(s)=x(t+s)$, $x_{t}\in \mathcal{C}$. Let $ X$ be the
closure in $\mathcal{C}$ of the set $\{f_{t},t\in
\mathbf{R}_{+}\}$, where $f:\mathbf{R}_{+}\rightarrow
\mathbf{R}_{+}^{*}$ is a decreasing function. It follows that $(
X,d)$ is a metric space. The mapping $\varphi: T\times
X\rightarrow X, \ \varphi(t,s,x)=x_{t-s}$ is an evolution semiflow
on $X$.

We consider $V=\mathbf{R}^{2}$, with the norm $\left\Vert
v\right\Vert=|v_{1}|+|v_{2}|$, $v=(v_{1},v_{2})\in V$. The mapping
$\Phi: T\times X\rightarrow \mathcal{B}(V)$ given by
\[
\Phi(t,s,x)v=\left(
e^{\alpha_{1}\int_{s}^{t}x(\tau-s)d\tau}v_{1},e^{\alpha_{2}\int_{s}^{t}x(\tau-s)d\tau}v_{2}\right),
\ (\alpha_{1},\alpha_{2})\in\mathbf{R}^{2},
\]
is an evolution cocycle over $\varphi$ and $C=(\varphi,\Phi)$ is a
skew-evolution semiflow.
\end{example}

\begin{remark}\rm
A connection between the solutions of a differential equation
\begin{equation}\label{ec_nuet1}
\dot{u}(t)=A(t)u(t), \ t\in\mathbf{R}_{+}
\end{equation}
and a skew-evolution semiflow is given by the definition of the
evolution cocycle $\Phi$, by the relation $\Phi(t,s,x)v=U(t,s)v$,
where $U(t,s)=u(t)u^{-1}(s)$, $(t,s)\in T$, $(x,v)\in Y$, and
where $u(t)$, $t\in \mathbf{R}_{+}$, is a solution of the
differential equation (\ref{ec_nuet1}).
\end{remark}

The fact that the skew-evolution semiflows are generalizations for
skew-product semiflows is emphasized by

\begin{example}\rm Let $X$ be
the metric space defined as in Example \ref{ex_ses}. The mapping
$\varphi_{0}:\mathbf{R}_{+}\times X\rightarrow X$,
$\varphi_{0}(t,x)=x_{t}$, where $x_{t}(\tau)=x(t+\tau)$, $\forall
\tau\geq 0$, is a semiflow on $X$. Let us consider for every $x\in
X$ the parabolic system with Neumann's boundary conditions:
\begin{equation}\label{sistem_parabolic}
\left\{
\begin{array}{lc}
\displaystyle \frac{\partial v}{\partial
t}(t,y)=x(t)\frac{\partial^{2}v}{\partial y^{2}}(t,y), & t>0, y\in (0,1) \\
v(0,y)=v_{0}(y), & y\in (0,1) \\
\displaystyle\frac{\partial v}{\partial y}(t,0)=\frac{\partial
v}{\partial
y}(t,1)=0, & t>0.%
\end{array}%
\right.
\end{equation}
Let $V=\mathcal{L}^{2}(0,1)$ be a separable Hilbert space with the
orthonormal basis
 $\{e_{n}\}_{n\in \mathbf{N}}$, $e_{0}=1$, $e_{n}(y)=\sqrt{2}\cos n\pi y$,
 where $y\in (0,1)$, $n\in \mathbf{N}$. We denote
$D(A)=\{v\in \mathcal{L}^{2}(0,1), \ v(0)=v(1)=0\}$ and we define
the operator
\[
A:D(A)\subset V\rightarrow V, \ Av=\frac{d^{2}v}{dy^{2}},
\]
which generates a $\mathcal{C}_{0}$-semigroup $S$, defined by
$S(t)v=\sum\limits_{n=0}^{\infty}e^{-n^{2}\pi^{2}t}\langle
v,e_{n}\rangle e_{n}$, where $\langle \cdot ,\cdot \rangle $
denotes the scalar product in $V$. We define for every $x\in X$,
$A(x):D(A)\subset V\rightarrow %
V$, $A(x)=x(0)A$, which allows us to rewrite system
(\ref{sistem_parabolic}) in $V$ as
\begin{equation}\label{sistem_modif}
\left\{
\begin{array}{lc}
\dot{v}(t)=A(\varphi_{0}(t,x))v(t), & t> 0 \\
v(0)=v_{0}.%
\end{array}%
\right.
\end{equation}
The mapping
\[
\Phi_{0}:\mathbf{R}_{+}\times X\rightarrow \mathcal{B}(V), \
\Phi_{0}(t,x)v=S\left(\int_{0}^{t}x(s)ds\right)v
\]
is a cocycle over the semiflow $\varphi_{0}$ and
$C_{0}=(\varphi_{0}, \Phi_{0})$ is a linear skew-product semiflow
strongly continuous on $Y$. Also, for all $v_{0}\in D(A)$, we have
obtained that $v(t)=\Phi (t,x)x_{0}, \ t\geq 0$, is a strongly
solution of system (\ref{sistem_modif}).

As $C_{0}=(\varphi_{0},\Phi_{0})$ is a skew-product semiflow on
$Y$, then the mapping $C: T\times Y\rightarrow Y$,
$C(t,s,x,v)=(\varphi(t,s,x),\Phi(t,s,x)v)$, where
\[
 \varphi(t,s,x)=\varphi_{0}(t-s,x) \ \textrm{and} \
 \Phi(t,s,x)=\Phi_{0}(t-s,x), \ \forall (t,s,x)\in  T\times X
\]
is a skew-evolution semiflow on $Y$. Hence, the skew-evolution
semiflows generalize the notion of skew-evolution semiflows.
\end{example}

An interesting class of skew-evolution semiflows, useful to
describe some asymptotic properties, is given by

\begin{example}\rm\label{ex_shift}
Let us consider a skew-evolution semiflow $C=(\varphi, \Phi)$ and
a parameter $\lambda \in \mathbf{R}$. We define the mapping
\begin{equation}\label{relcevshift}
\Phi_{\lambda}: T\times  X\rightarrow \mathcal{B}(V), \
\Phi_{\lambda}(t,t_{0},x)=e^{\lambda(t-t_{0})}\Phi(t,t_{0},x).
\end{equation}
One can remark that $C_{\lambda}=(\varphi, \Phi_{\lambda})$ also
satisfies the conditions of Definition \ref{def_coc_ev_1}, being
called \emph{$\lambda$-shifted skew-evolution semiflow} on $Y$.

Let us consider on the Banach space $V$ the Cauchy problem
\[
\left\{
\begin{array}{l}
\dot{v}(t)=Av(t), \ t> 0 \\
v(0)=v_{0}%
\end{array}%
\right.
\]
with the nonlinear operator $A$. Let us suppose that $A$ generates
a nonlinear $C_{0}$-semigroup $\mathcal{S}=\{S(t)\}_{t\geq 0}$.
Then $\Phi(t,s,x)v=S(t-s)v$, where $t\geq s\geq 0$, $(x,v)\in Y$,
defines an evolution cocycle. Moreover, the mapping defined by
$\Phi_{\lambda}: T\times  X\rightarrow \mathcal{B}(V)$,
$\Phi_{\lambda}(t,s,x)v=S_{\lambda}(t-s)v$, where
$\mathcal{S}_{\lambda}=\{S_{\lambda}(t)\}_{t\geq 0}$ is generated
by the operator $A-\lambda I$, is also an evolution cocycle.
\end{example}

\begin{definition}\rm\label{def_taremas}
A skew-evolution semiflow $C =(\varphi,\Phi)$ is said to be
\emph{strongly measurable} if, for all $(t_{0},x,v)\in T\times Y$,
the mapping $s\mapsto\left\Vert\Phi(s,t_{0},x)v\right\Vert$ is
measurable on $[t_{0},\infty)$.
\end{definition}

\begin{definition}\rm\label{def_neg}
The skew-evolution semiflow $C=(\varphi,\Phi)$ is said to have
\emph{exponential growth} if there exist
$M,\omega:\mathbf{R}_{+}\rightarrow\mathbf{R}_{+}^{*}$ such that:
\[
\left\Vert \Phi(t,t_{0},x)v\right\Vert \leq M(s)e^{\omega
(t-s)}\left\Vert \Phi(s,t_{0},x)v\right\Vert, \forall
(t,s),(s,t_{0})\in  T, \forall (x,v)\in Y.
\]%
\end{definition}

\begin{remark}\rm\label{obs_shift}
If $C=(\varphi,\Phi)$ is a skew-evolution semiflow with
exponential growth, as following relations
\[
\left\Vert\Phi_{\lambda}(t,t_{0},x)v\right\Vert=e^{\lambda(t-t_{0})}
\left\Vert\Phi(t,t_{0},x)v\right\Vert \leq
M(t_{0})e^{[\omega(t_{0})+\lambda](t-t_{0})}\left\Vert
v\right\Vert,
\]
hold for all $(t_{0},x,v)\in\mathbf{R}_{+}\times Y$, then
$C_{\lambda}=(\varphi,\Phi_{\lambda})$, $\lambda >0$, has also
exponential growth.
\end{remark}

\begin{remark}\rm\label{obs_eg}
$(i)$ If we consider in Definition \ref{def_neg} the constants
$M\geq 1$ and $\omega>0$, the skew-evolution semiflow $C$ is said
to have \emph{uniform exponential growth};

$(ii)$ If in Definition \ref{def_neg} we consider $M\geq 1$ to be
a constant such that the relation $\left\Vert
\Phi(t,s,x)\right\Vert \leq Me^{\omega (t-s)}$ holds for all
$(t,s)\in  T$ and all $x\in X$, the skew-evolution semiflow $C$ is
said to have \emph{bounded exponential growth}.
\end{remark}

\begin{definition}\rm\label{def_nedc}
The skew-evolution semiflow $C=(\varphi,\Phi)$ is said to have
\emph{exponential decay} if there exist
$M,\omega:\mathbf{R}_{+}\rightarrow\mathbf{R}_{+}^{*}$ such that:
\[
\left\Vert \Phi(s,t_{0},x)v\right\Vert \leq M(t)e^{\omega
(t-s)}\left\Vert \Phi(t,t_{0},x)v\right\Vert, \forall
(t,s),(s,t_{0})\in  T, \forall (x,v)\in Y.
\]%
\end{definition}

\begin{remark}\rm
If $C=(\varphi,\Phi)$ be a skew-evolution semiflow with
exponential decay, as following relations
\[
\left\Vert\Phi_{-\lambda}(s,t_{0},x)v\right\Vert=e^{-\lambda(s-t_{0})}
\left\Vert\Phi(s,t_{0},x)v\right\Vert \leq
M(t)e^{[\omega(t)+\lambda](t-s)}\left\Vert
\Phi_{-\lambda}(t,t_{0},x)v\right\Vert,
\]
hold for all $(t,s),(s,t_{0})\in T$ and all $(x,v)\in Y$, then
$C_{-\lambda}=(\varphi,\Phi_{-\lambda})$, $\lambda>0$, has also
exponential decay.
\end{remark}

\begin{remark}\rm
If in Definition \ref{def_nedc} we consider $M\geq 1$ and
$\omega>0$ to be constants, the skew-evolution semiflow $C$ is
said to have \emph{uniform exponential decay}.
\end{remark}

\section{On various classes of dichotomy}

Let $C: T\times Y\rightarrow Y$,
$C(t,s,x,v)=(\varphi(t,s,x),\Phi(t,s,x)v)$ be a skew-evolution
semiflow on $Y$.

\begin{definition}\rm\label{proiector}
A continuous mapping $P:Y\rightarrow Y$ defined by:
\begin{equation}
P(x,v)=(x,P(x)v), \ \forall (x,v)\in Y,
\end{equation}
where $P(x)$ is a linear projection on $Y_{x}$, is called
\emph{projector} on $Y$.
\end{definition}

\begin{remark}\rm
The mapping $P(x):Y_{x}\rightarrow Y_{x}$ is linear and bounded
and satisfies the relation $P(x)P(x)=P^{2}(x)=P(x)$ for all $x\in
X.$
\end{remark}

For all projectors $P:Y\rightarrow Y$ we define the sets
\[
Im P=\{(x,v)\in Y,P(x)v=v\} \ \textrm{and} \ KerP=\{(x,v)\in
Y,P(x)v=0\}.
\]

\begin{remark}\rm
Let $P$ be a projector on $Y$. Then $ImP$ and $KerP$ are closed
subsets of $Y$ and for all $x\in
 X$ we have
\[
ImP(x)+KerP(x)=Y_{x} \ \textrm{and} \ ImP(x)\cap KerP(x)=\{0\}.
\]
\end{remark}

\begin{remark}\rm
If $P$ is a projector on $Y$, then the mapping
\begin{equation}
Q:Y\rightarrow Y, \ Q(x,v)=(x,v-P(x)v)
\end{equation}
is also a projector on $Y$, called \emph{the complementary
projector} of $P$.
\end{remark}

\begin{definition}\rm\label{prinv}
A projector $P$ on $Y$ is called \emph{invariant} relative to a
skew-evolution semiflow $C=(\varphi,\Phi)$ if following relation
holds:
\begin{equation}
P(\varphi(t,s,x))\Phi(t,s,x)=\Phi(t,s,x)P(x),
\end{equation}
for all $(t,s)\in  T$ and all $x\in  X$.
\end{definition}

\begin{remark}\rm
If the projector $P$ is invariant relative to a skew-evolution
semiflow $C$, then its complementary projector $Q$ is also
invariant relative to $C$.
\end{remark}

\begin{definition}\rm\label{comp_pr_dich}
A projector $P_{1}$ and its complementary projector $P_{2}$ are
said to be \emph{compatible} with a skew-evolution semiflow
$C=(\varphi,\Phi)$ if

$(d_{1})$ the projectors $P_{1}$ and $P_{2}$ are invariant on $Y$;

$(d_{2})$ for all $x\in  X$, the projections $P_{1}(x)$ and
$P_{2}(x)$ commute and the relation $P_{1}(x)P_{2}(x)=0$ holds.
\end{definition}

In what follows we will denote
\[
\Phi_{k}(t,t_{0},x)=\Phi(t,t_{0},x)P_{k}(x), \ \forall
(t,t_{0})\in  T, \ \forall x\in  X, \ \forall k\in \{1,2\}.
\]
We remark that $\Phi_{k}$, $k\in \{1,2\}$ are evolution cocycles
and
\[
C_{k}(t,s,x,v)=(\varphi(t,s,x),\Phi_{k}(t,s,x)v), \ \forall
(t,t_{0},x,v)\in  T\times Y, \ \forall k\in \{1,2\},
\]
are skew-evolution semiflows, over all evolution semiflows
$\varphi$ on $X$.

\begin{definition}\rm\label{def_ued}
The skew-evolution semiflow $C=(\varphi,\Phi)$ is called
\emph{uniformly exponentially dichotomic} if there exist two
projectors $P_{1}$ and $P_{2}$ compatible with $C$, some constants
$N_{1}\geq 1$, $N_{2}\geq 1$ and $\nu_{1}$, $\nu_{2}>0$ such that:
\begin{equation}\label{dich_stab}
e^{\nu_{1}(t-s)}\left\Vert \Phi_{1}(t,t_{0},x)v\right\Vert \leq
N_{1}\left\Vert \Phi_{1}(s,t_{0},x)v\right\Vert;
\end{equation}
\begin{equation}\label{dich_instab}
e^{\nu_{2}(t-s)}\left\Vert \Phi_{2}(s,t_{0},x)(x)v\right\Vert \leq
N_{2}\left\Vert \Phi_{2}(t,t_{0},x)(x)v\right\Vert,
\end{equation}
for all $(t,s),(s,t_{0})\in T$ and all $(x,v)\in Y$.
\end{definition}

\begin{remark}\rm
Without any loss of generality we can consider
\[
N=\max\{N_{1},N_{2}\} \ \textrm{and} \ \nu=\min
\{\nu_{1},\nu_{2}\}.
\]
We will call $N_{1}$, $N_{2}$, $\nu_{1}$, $\nu_{2}$, respectively
$N$, $\nu$ \emph{dichotomic characteristics}.
\end{remark}

In what follows we will define generalizations for skew-evolution
semiflows of some asymptotic properties given by L. Barreira and
C. Valls for evolution equations in \cite{BaVa_LNM}.

\begin{definition}\rm\label{def_BVed}
The skew-evolution semiflow $C=(\varphi,\Phi)$ is called
\emph{Barreira-Valls exponentially dichotomic} if there exist two
projectors $P_{1}$ and $P_{2}$ compatible with $C$, some constants
$N\geq 1$, $\alpha_{1}$, $\alpha_{2}>0$ and $\beta_{1}$,
$\beta_{2}>0$ such that:
\begin{equation}\label{BV_dich_stab}
\left\Vert \Phi_{1}(t,t_{0},x)v\right\Vert \leq Ne^{-\alpha_{1}
t}e^{\beta_{1} s}\left\Vert \Phi_{1}(s,t_{0},x)v\right\Vert;
\end{equation}
\begin{equation}\label{BV_dich_instab}
\left\Vert \Phi_{2}(s,t_{0},x)v\right\Vert \leq Ne^{-\alpha_{2}
t}e^{\beta_{2} s}\left\Vert \Phi_{2}(t,t_{0},x)v\right\Vert,
\end{equation}
for all $(t,s),(s,t_{0})\in T$ and all $(x,v)\in Y$.
\end{definition}

\begin{definition}\rm\label{def_ed}
The skew-evolution semiflow $C=(\varphi,\Phi)$ is called
\emph{exponentially dichotomic} if there exist two projectors
$P_{1}$ and $P_{2}$ compatible with $C$, some mappings $N_{1}$,
$N_{2}:\mathbf{R}_{+}\rightarrow \mathbf{R}_{+}^{\ast }$ and some
constants $\nu_{1}$, $\nu_{2}>0$ such that:
\begin{equation}
\left\Vert \Phi_{1}(t,t_{0},x)v\right\Vert \leq
N_{1}(s)e^{-\nu_{1}t}\left\Vert \Phi_{1}(s,t_{0},x)v\right\Vert;
\end{equation}
\begin{equation}
\left\Vert \Phi_{2}(s,t_{0},x)v\right\Vert \leq
N_{2}(s)e^{-\nu_{2}t}\left\Vert \Phi_{2}(t,t_{0},x)v\right\Vert,
\end{equation}
for all $(t,s),(s,t_{0})\in T$ and all $(x,v)\in Y$.
\end{definition}

Some immediate connections concerning the previously defined
asymptotic properties for skew-evolution semiflows are given by

\begin{remark}\rm\label{obs_ued_BVed}
$(i)$ A uniformly exponentially dichotomic skew-evolution semiflow
is Barreira-Valls exponentially dichotomic;

$(ii)$ Barreira-Valls exponentially dichotomic skew-evolution
semiflow is exponentially dichotomic.
\end{remark}

The reciprocal statements are not true, as shown in what follows.
Hence, the next example emphasizes a skew-evolution semiflow which
is Barreira-Valls exponentially dichotomic, but is not uniformly
exponentially dichotomic.

\begin{example}\rm\label{ex_BVed}
Let $f:\mathbf{R}_{+}\rightarrow(0,\infty)$ be a decreasing
function with the property that there exists
$\lim\limits_{t\rightarrow\infty}f(t)=l>0$. We will consider
$\lambda>f(0)$. Let
$\mathcal{C}=\mathcal{C}(\mathbf{R},\mathbf{R})$ be the metric
space of all continuous functions $x:\mathbf{R}\rightarrow
\mathbf{R}$, with the topology of uniform convergence on compact
subsets of $\mathbf{R}$. $\mathcal{C}$ is metrizable relative to
the metric given in Example \ref{ex_ses}. We denote $ X$ the
closure in $\mathcal{C}$ of the set ${\{f_{t}, \ t\in
\mathbf{R}_{+}\}}$, where $f_{t}(\tau)=f(t+\tau)$, $\forall
\tau\in \mathbf{R}_{+}$. Then $( X,d)$ is a metric space. The
mapping $$\varphi: T\times X\rightarrow X, \
\varphi(t,s,x)(\tau)=x_{t-s}(\tau)=x(t-s+\tau)$$ is an evolution
semiflow on $ X$. Let us consider the Banach space
$V=\mathbf{R}^{2}$ with the norm $\left\Vert
v\right\Vert=|v_{1}|+|v_{2}|$, $v=(v_{1},v_{2})\in V$. The mapping
\[
\Phi: T\times  X \rightarrow \mathcal{B}(V), \ \Phi(t,s,x)v=
\]
\[
=\left(e^{t\sin t-s\sin
s-2(t-s)-\int_{s}^{t}x(\tau-s)d\tau}v_{1},\ e^{3(t-s)-2t\cos
t+2s\cos s+\int_{s}^{t}x(\tau-s)d\tau}v_{2}\right),
\]
where $t\geq s\geq 0, \ (x,v)\in Y$, is an evolution cocycle over
the evolution semiflow $\varphi$. We consider the projectors
$P_{1}, P_{2}:Y\rightarrow Y$, $P_{1}(x,v)=(v_{1},0)$, $
P_{2}(x,v)=(0,v_{2})$, for all $x\in X$ and all
$v=(v_{1},v_{2})\in V$, compatible with the skew-evolution
semiflow $C=(\varphi, \Phi)$.

We have, according to the properties of function $x$,
$$\left| \Phi(t,s,x)P_{1}(x)v\right|= e^{t\sin t-s\sin
s+2s-2t}e^{-\int_{s}^{t}x(\tau-s)d\tau}|v_{1}|\leq$$ $$\leq
e^{-t+3s}e^{-l(t-s)}|v_{1}|=e^{-(1+l)t}e^{(3+l)s}|v_{1}|,$$ for
all $(t,s,x,v)\in T\times Y$.

Also, following relations $$\left|
\Phi(t,s,x)P_{2}(x)v\right|=e^{3t-3s-2t\cos t+2s\cos
s+\int_{s}^{t}x(\tau-s)d\tau}|v_{2}|\geq$$ $$ \geq
e^{t-s}e^{l(t-s)}|v_{2}|=e^{(1+l)t}e^{-(1+l)s}|v_{2}|,$$ hold for
all $(t,s,x,v)\in T\times Y$.

Hence, the skew-evolution semiflow $C=(\varphi,\Phi)$ is
Barreira-Valls exponentially dichotomic with $N=1$,
$\alpha_{1}=\alpha_{2}=\beta_{2}=1+l$, $\beta_{1}=3+l$.

Let us suppose now that $C=(\varphi,\Phi)$ is uniformly
exponentially dichotomic. According to Definition \ref{def_ued},
there exist $N\geq 1$ and $\nu_{1}>0$, $\nu_{2}>0$ such that
$$e^{t\sin t-s\sin s+2s-2t}e^{-\int_{s}^{t}x(\tau-s)d\tau}|v_{1}|\leq Ne^{-\nu_{1}(t-s)}|v_{1}|, \ \forall t\geq s\geq 0$$
and $$Ne^{3t-3s-2t\cos t+2s\cos
s}e^{\int_{s}^{t}x(\tau-s)d\tau}|v_{2}|\geq
e^{\nu_{2}(t-s)}|v_{2}|, \ \forall t\geq s\geq 0.$$ If we consider
$t=2n\pi+\frac{\pi}{2}$ and $s=2n\pi$, we have in the first
inequality
$$e^{2n\pi-\frac{\pi}{2}}\leq Ne^{-\nu\frac{\pi}{2}}e^{\int\limits_{2n\pi}^{2n\pi+\frac{\pi}{2}}x(\tau-2n\pi)d\tau}
\leq Ne^{(-\nu_{1}+\lambda)\frac{\pi}{2}},$$ which, for
$n\rightarrow \infty$, leads to a contradiction. In the second
inequality, if we consider $t=2n\pi$ and $s=2n\pi-\pi$, we obtain
$$Ne^{-4n\pi+3\pi}\geq e^{\nu_{2}\pi}e^{\int\limits_{2n\pi-\pi}^{2n\pi}
x(\tau-2n\pi+\pi)d\tau}\geq e^{(\nu_{2}-\lambda)\pi}.$$ For
$n\rightarrow \infty$, a contradiction is obtained.

We obtain that $C$ is not uniformly exponentially dichotomic.
\end{example}

There exist exponentially dichotomic skew-evolution semiflows that
are not Barreira-Valls exponentially dichotomic, as in the next

\begin{example}\rm\label{ex_ed}
We consider the metric space $(X,d)$, the Banach space $V$, the
evolution semiflow $\varphi$ and the projectors $P_{1}$ and
$P_{2}$ defined as in Example \ref{ex_BVed}. Let us consider a
continuous function
\[
g:\mathbf{R}_{+}\rightarrow[1,\infty) \ \textrm{with} \
g(n)=e^{n\cdot2^{2n}} \ \textrm{and} \
g\left(n+\frac{1}{2^{2n}}\right)=1.
\]
The mapping $\Phi: T\times X\rightarrow \mathcal{B}(V)$, defined
by
\[
\Phi(t,s,x)v=\left(\frac{g(s)}{g(t)}e^{-(t-s)-\int_{s}^{t}x(\tau-s)d\tau}v_{1},
\frac{g(t)}{g(s)}e^{t-s+\int_{s}^{t}x(\tau-s)d\tau}v_{2}\right)
\]
is an evolution cocycle over the evolution semiflow $\varphi$. As
$$\left| \Phi(t,s,x)P_{1}(x)v\right|\leq g(s)e^{-(1+l)(t-s)}|v_{1}|, \ \forall
(t,s,x,v)\in T\times Y$$ and $$g(s)\left|
\Phi(t,s,x)P_{2}(x)v\right|\geq e^{(1+l)(t-s)}|v_{2}|,\ \forall
(t,s,x,v)\in T\times Y,$$ the skew-evolution semiflow
$C=(\varphi,\Phi)$ is exponentially dichotomic, with
$N_{1}(u)=N_{2}(u)=g(u)\cdot e^{(1+l)u}$, $u\geq 0$, and
$\nu_{1}=\nu_{2}=1+l$.

Let us suppose that $C$ is Barreira-Valls exponentially
dichotomic. There exist $N\geq 1$ and $\alpha_{1}$, $\alpha_{2}$,
$\beta_{1}$, $\beta_{2}>0$ such that
$$\frac{g(s)}{g(t)}e^{s}\leq Ne^{t}e^{-\alpha_{1} t}e^{\beta_{1} s}e^{\int_{s}^{t}x(\tau-s)d\tau}$$
and $$e^{\alpha_{2} t}e^{-t}\leq N\frac{g(t)}{g(s)}e^{\beta_{2}
s}e^{-s}e^{\int_{s}^{t}x(\tau-s)d\tau}.$$ Further, if we consider
$t=n+\displaystyle\frac{1}{2^{2n}}$ and $s=n$, it follows that
$$e^{n(2^{2n}+1+\alpha_{1}-\beta_{1})}\leq Ne^{\frac{\lambda-\alpha_{1}}{2^{2n}}}\ \textrm{and}
\ e^{n(2^{n}+\alpha_{2}-\beta_{2})}\leq
Ne^{\frac{1+\lambda-\alpha_{2}}{2^{2n}}}.$$ As, for $n\rightarrow
\infty$, two contradictions are obtained, it follows that $C$ is
not Barreira-Valls exponentially dichotomic.
\end{example}

Let us present some particular classes of dichotomy, given by

\begin{definition}\rm\label{upd}
A skew-evolution semiflow $C =(\varphi,\Phi)$ is \emph{uniformly
polynomially dichotomic} if there exist two projectors $P_{1}$ and
$P_{2}$ compatible with $C$ and some constants $N\geq 1$ and
$\alpha_{1}>0$, $\alpha_{2}>0$ such that:
\begin{equation}
\left\Vert \Phi_{1}(t,s,x)v\right\Vert \leq Nt^{-\alpha_{1}
}s^{\alpha_{1}}\left\Vert P_{1}(x)v\right\Vert;
\end{equation}
\begin{equation}
\left\Vert P_{2}(x)v\right\Vert \leq Nt^{-\alpha_{2}
}s^{\alpha_{2}}\left\Vert \Phi_{2}(t,s,x)v\right\Vert;
\end{equation}
for all $(t,s)\in T$ and all $(x,v)\in Y$.
\end{definition}

\begin{definition}\rm\label{BVpd}
A skew-evolution semiflow $C =(\varphi,\Phi)$ is
\emph{Barreira-Valls polynomially dichotomic} if there exist some
constants $N\geq 1$, $\alpha_{1}>0$, $\alpha_{2}>0$ and
$\beta_{1}>0$, $\beta_{2}>0$ such that:
\begin{equation}
\left\Vert \Phi_{1}(t,s,x)v\right\Vert \leq Nt^{-\alpha_{1}
}s^{\beta_{1}}\left\Vert P_{1}(x)v\right\Vert;
\end{equation}
\begin{equation}
\left\Vert P_{2}(x)v\right\Vert \leq Nt^{-\alpha_{2}
}s^{\beta_{2}}\left\Vert \Phi_{2}(t,s,x)v\right\Vert,
\end{equation}
for all $(t,s)\in T$ and all $(x,v)\in Y$.
\end{definition}

\begin{definition}\rm\label{pd}
A skew-evolution semiflow $C =(\varphi,\Phi)$ is
\emph{polynomially dichotomic} if there exist a function
$N:\mathbf{R}_{+}\rightarrow [1,\infty)$, some constants
$\alpha_{1}>0$ and $\alpha_{2}>0$ such that:
\begin{equation}
\left\Vert \Phi_{1}(t,s,x)v\right\Vert \leq N(s)t^{-\alpha_{1}
}\left\Vert P_{1}(x)v\right\Vert;
\end{equation}
\begin{equation}
\left\Vert P_{2}(x)v\right\Vert \leq N(s)t^{-\alpha_{2}
}\left\Vert \Phi_{2}(t,s,x)v\right\Vert,
\end{equation}
for all $(t,s)\in T$ and all $(x,v)\in Y$.
\end{definition}

Relations between the defined classes of dichotomy are described
by

\begin{remark}\rm\label{obs_upd_BVpd}
$(i)$ A uniformly polynomially dichotomic skew-evolution semiflow
is Barreira-Valls polynomially dichotomic;

$(ii)$ A Barreira-Valls polynomially dichotomic is polynomially
dichotomic.
\end{remark}

The next example shows a skew-evolution semiflow which is
Barreira-Valls polynomially dichotomic but is not uniformly
polynomially dichotomic.

\begin{example}\rm
We consider the metric space $(X,d)$, the Banach space $V$, the
evolution semiflow $\varphi$ and the projectors $P_{1}$ and
$P_{2}$ defined as in Example \ref{ex_BVed}. We will consider the
mapping
$$g:\mathbf{R}_{+}\rightarrow\mathbf{R}, \ g(t)=(t+1)^{3-\sin\ln(t+1)}.$$ We define
$$\Phi(t,s,x)v=\left(\frac{g(s)}{g(t)}e^{-\int_{s}^{t}x(\tau-s)d\tau}v_{1},
\frac{g(t)}{g(s)}e^{\int_{s}^{t}x(\tau-s)d\tau}v_{2}\right),
(t,s)\in T,\ (x,v)\in Y.$$ $\Phi$ is an evolution cocycle over
$\varphi$. Due to the properties of function $x$ and of function
$f:(0,\infty)\rightarrow (0,\infty), \
f(u)=\displaystyle\frac{e^{u}}{u}$, we have
$$\left| \Phi_{1}(t,s,x)v\right| \leq \frac{(s+1)^{4}}{(t+1)^{2}}e^{-l(t-s)}|v_{1}|\leq (s+1)^{2}
\left(\frac{s+1}{t+1}\right)^{2}e^{-lt}e^{ls}|v_{1}|\leq$$ $$\leq
\frac{s(s+1)^{2}}{t}t^{-l}s^{l}|v_{1}|\leq
4t^{-(1+l)}s^{3+l}|v_{1}|,$$ for all $t\geq s\geq t_{0}=1$ and all
$(x,v)\in Y$. Also, following relations $$\left|
\Phi_{2}(t,s,x)v\right| \geq
\frac{(s+1)^{4}}{(t+1)^{2}}e^{-l(t-s)}|v_{2}|\geq
\frac{(t+1)^{2}}{(s+1)^{4}}e^{lt}e^{-ls}|v_{2}|\geq
t^{2+l}s^{-8-l}|v_{2}|,$$ hold for all $t\geq s\geq t_{0}=1$ and
all $(x,v)\in Y$.

Hence, by Definition \ref{BVpd}, the skew-evolution semiflow
$C=(\varphi,\Phi )$ is Barreira-Valls polynomially dichotomic.

We suppose now that $C$ is uniformly polynomially dichotomic.
According to Definition \ref{upd}, there exist $N\geq 1$ and
$\alpha_{1}>0$ such that
$$\frac{(s+1)^{3}}{(t+1)^{3}}\frac{(t+1)^{\sin\ln (t+1)}}{(s+1)^{\sin\ln (s+1)}}
\leq
Nt^{-\alpha_{1}}s^{\alpha_{1}}e^{\int_{s}^{t}x(\tau-s)d\tau}$$ for
all $t\geq s\geq t_{0}$. Let us consider
$$t=e^{2n\pi+\frac{\pi}{2}}-1 \ \textrm{and} \
s=e^{2n\pi-\frac{\pi}{2}}-1.$$ We have, if we consider the
properties of function $x$, that
$$e^{(2n-\lambda-1)\pi}\leq N e^{2\alpha_{1}},$$
which, if $n\rightarrow \infty$, leads to a contradiction.

Also, as in Definition \ref{upd}, there exist $N\geq 1$ and
$\alpha_{1}>0$ such that
$$N\frac{(t+1)^{3}}{(s+1)^{3}}\frac{(s+1)^{\sin\ln (s+1)}}{(t+1)^{\sin\ln (t+1)}}
\geq
t^{\alpha_{2}}s^{-\alpha_{2}}e^{-\int_{s}^{t}x(\tau-s)d\tau}$$ for
all $t\geq s\geq t_{0}$, which implies, for
$t=e^{2n\pi+\frac{\pi}{2}}-1$ and $s=e^{2n\pi-\frac{\pi}{2}}-1$,
$$Ne^{(-2n+\lambda-1)\pi}\geq e^{-2\alpha_{2}},$$ which, for $n\rightarrow \infty$, is a contradiction.

We obtain thus that $C$ is not uniformly polynomially dichotomic.
\end{example}

There exist skew-evolution semiflows that are polynomially
dichotomic but are not Barreira-Valls polynomially dichotomic.

\begin{example}\rm
Let us consider the data given in Example \ref{ex_ed}. We obtain
$$\left| \Phi(t,s,x)P_{1}(x)v\right|\leq g(s)e^{-(1+l)(t-s)}|v_{1}|
\leq g(s)e^{(1+l)s}t^{-(1+l)}$$ and $$g(s)e^{(1+l)s}\left|
\Phi(t,s,x)P_{2}(x)v\right|\geq e^{(1+l)t}|v_{2}|\geq
t^{(1+l)}|v_{2}|,$$ for all $(t,s,x,v)\in T\times Y$, which proves
that the skew-evolution semiflow $C=(\varphi,\Phi)$ is
polynomially dichotomic.

If we suppose that $C$ is Barreira-Valls polynomially dichotomic,
there exist $N\geq 1$, $\alpha_{1}>0$, $\alpha_{2}>0$ and
$\beta_{1}>0$, $\beta_{2}>0$ such that $$\frac{g(s)}{g(t)}\leq
Nt^{-\alpha_{1}}s^{\beta_{1}}e^{t-s+\int_{s}^{t}x(\tau-s)d\tau}\
\textrm{and} \ t^{\alpha_{2}}\leq N\frac{g(t)}{g(s)}s^{\beta_{2}
}e^{t-s+\int_{s}^{t}x(\tau-s)d\tau}.$$ If we consider
$t=n+\displaystyle\frac{1}{2^{2n}}$ and $s=n$, we obtain
$$e^{n\cdot 2^{2n}}\leq N\cdot n^{-\alpha_{1}}\cdot
n^{\beta_{1}}\cdot e^{\frac{1+\lambda}{2^{2n}}}\ \textrm{and}\
e^{n\cdot 2^{2n}} \leq
N\left(n+\frac{1}{2^{2n}}\right)^{-\alpha_{2}}\cdot
n^{\beta_{2}}\cdot e^{\frac{1+\lambda}{2^{2n}}}.$$ For
$n\rightarrow \infty$, two contradictions are obtained, which
proves that $C$ is not Barreira-Valls polynomially dichotomic.
\end{example}

\section{Main results}

The first results will prove some relations between all the
classes of dichotomies.

\begin{proposition}
A uniformly exponentially dichotomic skew-evolution semiflow
$C=(\varphi, \Phi)$ is uniformly polynomially dichotomic.
\end{proposition}

\begin{proof}
Let us consider in Definition \ref{def_ued}, without any loss of
generality, $t_{0}=1$. It also assures the existence of constants
$N\geq 1$ and $\nu_{1} >0$ such that $\left\Vert
\Phi_{1}(t,s,x)v\right\Vert \leq Ne^{-\nu_{1} (t-s)}\left\Vert
P_{1}(x)v\right\Vert.$ As
$$e^{-u}\leq \frac{1}{u+1}, \ \forall u\geq 0 \ \textrm{and }\
\frac{t}{s}\leq t-s+1, \ \forall t\geq s\geq 1,$$ it follows that
$$\left\Vert \Phi_{1}(t,s,x)v\right\Vert \leq
N(t-s+1)^{-\nu_{1}}\left\Vert P_{1}(x)v\right\Vert\leq
Nt^{-\nu_{1}}s^{\nu_{1}}\left\Vert P_{1}(x)v\right\Vert,$$ for all
$t\geq s\geq 1$ and all $(x,v)\in Y$.

We also have the property of function $$f:(0,\infty)\rightarrow
(0,\infty), \ f(u)=\frac{e^{u}}{u}$$ of being nondecreasing, which
assures the inequality $$\frac{e^{s}}{e^{t}}\leq\frac{s}{t},\
\forall t\geq s>0$$ and, further, for all $t\geq s\geq 1$ and all
$(x,v)\in Y$, we have
$$\left\Vert P_{2}(x)v\right\Vert\leq Ne^{-\nu_{2}t}e^{\nu_{2}
s}\left\Vert \Phi_{2}(t,s,x)v\right\Vert \leq
Nt^{-\nu_{2}}s^{\nu_{2} }\left\Vert \Phi_{2}(t,s,x)v\right\Vert,$$
where constants $N\geq 1$ and $\nu_{2} >0$ are also given by
Definition \ref{def_ued}.

Thus, according to Definition \ref{upd}, $C$ is uniformly
polynomially dichotomic.
\end{proof}

\vspace{3mm}

We give an example of a skew-evolution semiflow which is uniformly
polynomially dichotomic, but is not uniformly exponentially
dichotomic.

\begin{example}\rm\label{ex_upd}
Let $(X,d)$ be the metric space, $V$ the Banach space, $\varphi$
the evolution semiflow, $P_{1}$ and $P_{2}$ the projectors given
as in Example \ref{ex_BVed}.

Let us consider the function
$g:\mathbf{R}_{+}\rightarrow\mathbf{R}$, given by $g(t)=t^{2}+1$
and let us define
$$\Phi(t,s,x)v=\left(\frac{g(s)}{g(t)}e^{-\int_{s}^{t}x(\tau-s)d\tau}v_{1},
\frac{g(t)}{g(s)}e^{\int_{s}^{t}x(\tau-s)d\tau}v_{2}\right), \
(t,s)\in T,\ (x,v)\in Y.$$ We can consider $t_{0}=1$ in Definition
\ref{upd}. As,
$\displaystyle\frac{s^{2}+1}{t^{2}+1}\leq\frac{s}{t}$, for $t\geq
s\geq 1$ and according to the properties of function $x$, we have
$$\frac{s^{2}+1}{t^{2}+1}e^{-\int_{s}^{t}x(\tau-s)d\tau}|v_{1}|
\leq t^{-(1+l)}s^{1+l}|v_{1}|$$ and $$\frac{t^{2}+1}{s^{2}+1}
e^{\int_{s}^{t}x(\tau-s)d\tau}|v_{2}| \geq
t^{(2+l)}s^{-(4+l)}|v_{2}|,$$ for all $t\geq s\geq 1$ and all
$v\in V$. It follows that $C=(\varphi,\Phi)$ is uniformly
polynomially dichotomic.

If the skew-evolution semiflow $C=(\varphi,\Phi)$ is also
uniformly exponentially dichotomic, according to Definition
\ref{def_ued}, there exist $N\geq 1$ $\nu_{1}>0$ and $\nu_{2}>0$
such that
$$\frac{s^{2}+1}{t^{2}+1}|v_{1}|\leq Ne^{-\nu_{1}(t-s)}e^{-l(t-s)}|v_{1}|\
\textrm{and}\ N\frac{t^{2}+1}{s^{2}+1}|v_{2}|\geq
e^{\nu_{2}(t-s)}e^{l(t-s)}|v_{2}|,$$ for all $t\geq s\geq t_{0}$
and all $v\in V.$ If we consider $s=t_{0}$ and $t\rightarrow
\infty$, two contradictions are obtained, which proves that $C$ is
not uniformly exponentially dichotomic.
\end{example}

\begin{proposition}
A Barreira-Valls exponentially dichotomic skew-evolution semiflow
$C=(\varphi, \Phi)$ with $\alpha_{i}\geq\beta_{i}>0$, $i\in
\{1,2\}$, is Barreira-Valls polynomially dichotomic.
\end{proposition}

\begin{proof}
According to Definition \ref{def_BVed}, there exist some constants
$N\geq 1$, $\alpha_{1}>0$ and $\beta_{1}>0$ such that $$\left\Vert
\Phi_{1}(t,s,x)v\right\Vert \leq Ne^{-\alpha_{1} t}e^{\beta_{1}
s}\left\Vert P_{1}(x)v\right\Vert, \ \forall (t,s)\in T,\ \forall
(x,v)\in Y.$$ As the mapping $f:(0,\infty)\rightarrow (0,\infty)$,
defined by $f(u)=\displaystyle\frac{e^{u}}{u}$ is nondecreasing,
and as, by hypothesis, we can chose $\alpha_{1}\geq\beta_{1}$, we
obtain that
$$\left\Vert \Phi_{1}(t,s,x)v\right\Vert \leq
Nt^{-\alpha_{1}}e^{-\beta_{1} s}s^{\beta_{1}}e^{\beta_{1}
s}\left\Vert
P_{1}(x)v\right\Vert=Nt^{-\alpha_{1}}s^{\beta_{1}}\left\Vert
P_{1}(x)v\right\Vert,$$ for all $t\geq s>0$ and all $(x,v)\in Y$.

Analogously, we obtain $$\left\Vert P_{2}(x)v\right\Vert\leq
Ne^{-\alpha_{2}t}e^{\beta_{2}s}\left\Vert
\Phi_{2}(t,s,x)v\right\Vert\leq
Nt^{-\alpha_{2}}s^{\beta_{2}}\left\Vert
\Phi_{2}(t,s,x)v\right\Vert,$$ for all $t\geq s>0$ and all
$(x,v)\in Y$, where the constants $N\geq 1$, $\alpha_{2}>0$ and
$\beta_{2}>0$ are also assured by Definition \ref{def_ued}, with
the property $\alpha_{2}\geq\beta_{2}$.

Hence, according to Definition \ref{upd}, $C$ is Barreira-Valls
polynomially dichotomic.
\end{proof}

\vspace{3mm}

There exist skew-evolution semiflows that are Barreira-Valls
polynomially dichotomic, but are not Barreira-Valls exponentially
dichotomic.

\begin{example}\rm\label{ex_BVpd}
We consider the metric space $(X,d)$, the Banach space $V$, the
evolution semiflow $\varphi$ and the projectors $P_{1}$ and
$P_{2}$ defined as in Example \ref{ex_BVed}.

Let us consider the function
$g:\mathbf{R}_{+}\rightarrow\mathbf{R}$, given by $g(t)=t+1$ and
let us define an evolution cocycle $\Phi$ as in Example
\ref{ex_upd}. We obtain
$$\frac{s+1}{t+1}e^{-\int_{s}^{t}x(\tau-s)d\tau}|v_{1}|
\leq \frac{s^{2}}{t}e^{-l(t-s)}|v_{1}|\leq
t^{-1-l}s^{2+l}|v_{1}|$$ and
$$\frac{t+1}{s+1}e^{\int_{s}^{t}x(\tau-s)d\tau}|v_{2}| \geq t^{1+l}s^{-2-l}|v_{2}|,$$
for all $t\geq s\geq 1$ and all $v\in V$. It follows that the
skew-evolution semiflow $C=(\varphi,\Phi)$ is Barreira-Valls
polynomially dichotomic.

Let us suppose that $C$ is also Barreira-Valls exponentially
dichotomic. According to Definition \ref{def_BVed}, there exist
some constants $N\geq 1$, $\alpha_{1},\beta_{1}>0$ and
$\alpha_{2},\beta_{2}>0$ such that
$$\frac{s+1}{t+1}e^{-\int_{s}^{t}x(\tau-s)d\tau}|v_{1}|
\leq Ne^{-\alpha_{1} t}e^{\beta_{1} s}|v_{1}|$$ and
$$N\frac{t+1}{s+1}e^{\int_{s}^{t}x(\tau-s)d\tau}|v_{2}| \geq
e^{\alpha_{2} t}e^{-\beta_{2} s}|v_{2}|,$$ for all
$(t,s),(s,t_{0})\in T$ and all $(x,v)\in Y$. We consider
$s=t_{0}$. We have
$$\frac{e^{\alpha_{1} t}}{t+1}\leq\frac{\overline{N}}{t_{0}+1}\ \textrm{and}
\ \frac{e^{\alpha_{2} t}}{t+1}\leq\frac{\widetilde{N}}{t_{0}+1}, \
\forall t\geq t_{0}.$$ For $t\rightarrow \infty$, we obtain two
contradictions, and, hence, $C$ is not Barreira-Valls
exponentially dichotomic.
\end{example}

\begin{proposition}
An exponentially dichotomic skew-evolution semiflow $C=(\varphi,
\Phi)$ is polynomially dichotomic.
\end{proposition}

\begin{proof}
Definition \ref{def_ed} assures the existence of a function
$N_{1}:\mathbf{R}_{+}\rightarrow[1,\infty)$ and a constant
$\nu_{1}>0$ such that $$\left\Vert \Phi_{1}(t,s,x)v\right\Vert\leq
N_{1}(s)e^{-\nu_{1} t}\left\Vert P_{1}(x)v\right\Vert, \ \forall
(t,s)\in T, \ \forall (x,v)\in Y.$$ As following inequalities
$e^{t}\geq t+1>t$ hold for all $t\geq 0$, we obtain $$\left\Vert
\Phi_{1}(t,s,x)v\right\Vert\leq N_{1}(s)t^{-\nu_{1}}\left\Vert
P_{1}(x)v\right\Vert,$$ for all $t\geq s>0$ and all $(x,v)\in Y$.

As, by Definition \ref{def_ed} there exist a function
$N_{2}:\mathbf{R}_{+}\rightarrow[1,\infty)$ and a constant
$\nu_{2}>0$ such that $$\left\Vert P_{2}(x)v\right\Vert\leq
N_{2}(s)e^{-\nu_{2} t}\left\Vert \Phi_{2}(t,s,x)v\right\Vert, \
\forall (t,s)\in T, \ \forall (x,v)\in Y.$$ Analogously, as
previously, we have $$\left\Vert P_{2}(x)v\right\Vert\leq
N_{2}(s)t^{-\nu_{2}}\left\Vert \Phi_{2}(t,s,x)v\right\Vert,$$ for
all $t\geq s>0$ and all $(x,v)\in Y$.

Hence, according to Definition \ref{pd}, $C$ is polynomially
dichotomic.
\end{proof}

\vspace{3mm}

We present an example of a skew-evolution semiflow which is
polynomially dichotomic, but is not exponentially dichotomic.

\begin{example}\rm
We consider the metric space $(X,d)$, the Banach space $V$, the
evolution semiflow $\varphi$, the projectors $P_{1}$, $P_{2}$ and
function $g$ as in Example \ref{ex_BVpd}. Let
$$\Phi(t,s,x)v=\left(\frac{g(s)}{g(t)}e^{\int_{s}^{t}x(\tau-s)d\tau}|v_{1}|,
\frac{g(t)}{g(s)}e^{-\int_{s}^{t}x(\tau-s)d\tau}|v_{2}|\right)$$
be an evolution cocycle. Analogously as in the mentioned Example,
the skew-evolution semiflow $C$ is Barreira-Valls polynomially
dichotomic, and, according to Remark \ref{obs_upd_BVpd} $(ii)$, it
is also polynomially dichotomic. On the other hand, if we suppose
that $C$ is exponentially dichotomic, there exist $N_{1}$,
$N_{2}:\mathbf{R}_{+}\rightarrow \mathbf{R}_{+}^{\ast }$ and
$\nu_{1}$, $\nu_{2}>0$ such that
$$\frac{s+1}{t+1}|v_{1}| \leq
N_{1}(s)e^{-(\nu_{1}+l)t}e^{l s}|v_{1}|$$ and $$ |v_{2}| \leq
N_{2}(s)e^{-(\nu_{2}+l)t} \frac{t+1}{s+1}|v_{2}|,$$ for all
$(t,s)\in T$ and all $(x,v)\in Y$. If we consider $s=t_{0}$ and
$t\rightarrow \infty$, we obtain two contradictions, which shows
that $C$ is not exponentially dichotomic.
\end{example}

A characterization for the classic and mostly encountered property
of exponential dichotomy is given by the next

\begin{theorem}
Let $C=(\varphi,\Phi)$ be a strongly measurable skew-evolution
semiflow. $C$ is exponentially dichotomic if and only if there
exist two projectors $P_{1}$ and $P_{2}$ compatible with $C$ with
the properties that $C_{1}$ has bounded exponential growth and
$C_{2}$ has exponential decay such that

$(i)$ there exist a constant $\gamma>0$ and a mapping
$D:\mathbf{R}_{+}\rightarrow [1,\infty)$ with the property:
$$\int_{s}^{\infty}e^{(\tau-s)\gamma}\left\Vert
\Phi_{1}(\tau,s,x)v\right\Vert d\tau \leq D(s)\left\Vert
P_{1}(x)v\right\Vert,$$ for all $s\geq 0$ and all $(x,v)\in Y$;

$(ii)$ there exist a constant $\rho>0$ and a nondecreasing mapping
$\widetilde{D}:\mathbf{R}_{+}\rightarrow [1,\infty)$ with the
property:
$$\int_{t_{0}}^{t}e^{(t-\tau)\rho}\left\Vert
\Phi_{2}(\tau,t_{0},x)v\right\Vert d\tau \leq
\widetilde{D}(t_{0})\left\Vert \Phi_{2}(t,t_{0},x)v\right\Vert,$$
for all $t\geq t_{0}\geq 0$ and all $(x,v)\in Y$.
\end{theorem}

\begin{proof}
\emph{Necessity.} As $C$ is exponentially dichotomic, according to
Definition \ref{def_ed}, there exist $N_{1}$,
$N_{2}:\mathbf{R}_{+}\rightarrow \mathbf{R}_{+}^{\ast }$ and
$\nu_{1}$, $\nu_{2}>0$ such that
$$\left\Vert \Phi_{1}(t,t_{0},x)v\right\Vert \leq
N_{1}(s)e^{-\nu_{1}t}\left\Vert \Phi_{1}(s,t_{0},x)v\right\Vert$$
and $$ \left\Vert \Phi_{2}(s,t_{0},x)v\right\Vert \leq
N_{2}(s)e^{-\nu_{2}t}\left\Vert \Phi_{2}(t,t_{0},x)v\right\Vert,$$
for all $(t,s),(s,t_{0})\in T$ and all $(x,v)\in Y$.

In order to prove $(i)$, let us define
$\gamma=-\displaystyle\frac{\nu_{1}}{2}$. We obtain successively
$$\int_{s}^{\infty}e^{(\tau-s)\gamma}\left\Vert
\Phi_{1}(\tau,s,x)v\right\Vert d\tau \leq N_{1}(s)\left\Vert
\Phi_{1}(s,s,x)v\right\Vert\int_{s}^{\infty}e^{-\frac{\nu_{1}}{2}(\tau-s)}e^{-\nu_{1}(s-\tau)}d\tau=$$
$$=N(s)\left\Vert
P_{1}(x)v\right\Vert\int_{s}^{\infty}e^{-\frac{\nu_{1}}{2}(s-\tau)}d\tau=D(s)\left\Vert
P_{1}(x)v\right\Vert,$$ for all $s\geq 0$ and all $(x,v)\in Y$,
where we have denoted $$D(u)=\frac{N_{1}(u)}{\gamma}, \ u\geq 0.$$

To prove $(ii)$, we define $\rho=\displaystyle\frac{\nu_{2}}{2}$.
Following relations $$\int_{t_{0}}^{t}e^{(t-\tau)\rho}\left\Vert
\Phi_{2}(\tau,t_{0},x)v\right\Vert d\tau \leq
N_{2}(t_{0})\left\Vert\Phi_{2}(t,t_{0},x)v\right\Vert
\int_{t_{0}}^{t}e^{\frac{\nu_{2}}{2}(t-\tau)}e^{-\nu_{2}(t-\tau)}d\tau\leq$$
$$\leq \widetilde{D}(t_{0})\left\Vert\Phi_{2}(t,t_{0},x)v\right\Vert$$
hold for all $s\geq 0$ and all $(x,v)\in Y$, where we have
denoted
$$\widetilde{D}(u)=\frac{2N_{2}(u)}{\rho}, \ u\geq 0.$$
\emph{Sufficiency.} According to relation $(i)$, the
$\gamma$-shifted skew-evolution semiflow
$C_{\gamma}^{1}=(\varphi,\Phi^{1}_{\gamma})$, defined as in
Example \ref{ex_shift}, has bounded exponential growth and there
exists $D:\mathbf{R}_{+}\rightarrow [1,\infty)$ such that
$$\int_{s}^{\infty}\left\Vert
\Phi^{1}_{\gamma}(\tau,s,x)v\right\Vert d\tau \leq D(s)\left\Vert
P_{1}(x)v\right\Vert,$$ for all $s\geq 0$ and all $(x,v)\in Y$.

First of all, we will prove that there exists
$D_{1}:\mathbf{R}_{+}\rightarrow [1,\infty)$ such that $\left\Vert
\Phi^{1}_{\gamma}(t,s,x)v\right\Vert \leq D_{1}(s)\left\Vert
P_{1}(x)v\right\Vert$, for all $t\geq s\geq 0$ and all $(x,v)\in
Y$. Let us consider, for $t\geq s+1$,
$$c=\int_{0}^{1}e^{-\omega(u)}d u\leq \int^{t-s}_{0}e^{-\omega(u)}d
u=\int^{t}_{s}e^{-\omega(t-\tau)}d\tau.$$ Hence, for $t\geq s+1$,
we obtain $$c|<v^{*},\Phi^{1}_{\gamma}(t,s,x)v>|\leq
\int^{t}_{s}e^{-\omega(t-\tau)}|<v^{*},\Phi^{1}_{\gamma}(t,s,x)v>|d\tau=$$
$$=\int^{t}_{s}e^{-\omega(t-\tau)}\left\Vert
\Phi^{1}_{\gamma}(t,\tau,\varphi(\tau,s,x))^{*}v^{*}\right\Vert
\left\Vert \Phi^{1}_{\gamma}(\tau,s,x)v\right\Vert d\tau\leq$$
$$\leq M\left\Vert
P_{1}(x)v^{*}\right\Vert\int_{s}^{t}\left\Vert
\Phi^{1}_{\gamma}(\tau,s,x)v\right\Vert d\tau\leq MD(s)\left\Vert
P_{1}(x)v\right\Vert\left\Vert P_{1}(x)v^{*}\right\Vert,$$ where
$v\in V$, $v^{*}\in V^{*}$ and $M$, $\omega$ are given by
Definition \ref{def_neg} and Remark \ref{obs_eg}. Hence,
$$\left\Vert \Phi_{1}(t,s,x)v\right\Vert\leq \frac{MD(s)}{c},
\ \forall t\geq s+1, \ \forall (x,v)\in Y.$$

Now, for $t\in[s,s+1)$, we have
$$\left\Vert \Phi_{1}(t,s,x)v\right\Vert\leq Me^{\omega(1)}\left\Vert P_{1}(x)v\right\Vert,
 \ \forall (x,v)\in Y.$$
Thus, we obtain $$\left\Vert \Phi^{1}_{\gamma}(t,s,x)v\right\Vert
\leq D_{1}(s)\left\Vert P_{1}(x)v\right\Vert,$$ for all $t\geq
s\geq 0$ and all $(x,v)\in Y$, where we have denoted $$D_{1}(u)=
M\left[e^{\omega(1)}+\frac{D(u)}{c}\right], \ u\geq 0.$$ Further,
it follows that
$$\left\Vert \Phi_{1}(t,s,x)v\right\Vert \leq
D_{1}(s)e^{-(t-s)\gamma}\left\Vert v\right\Vert,\ \forall t\geq
s\geq 0.$$

According to $(ii)$, there exist a constant $\rho>0$ and a
nondecreasing mapping $\widetilde{D}:\mathbf{R}_{+}\rightarrow
[1,\infty)$ such that
$$\int_{t_{0}}^{t}e^{-(\tau-t_{0})\rho}\left\Vert
\Phi_{2}(\tau,t_{0},x)v\right\Vert d\tau \leq
\widetilde{D}(t_{0})e^{-(t-t_{0})\rho}\left\Vert
\Phi_{2}(t,t_{0},x)v\right\Vert,$$ for all $t\geq t_{0}\geq 0$ and
all $(x,v)\in Y$. Thus,
$$\int_{t_{0}}^{t}\left\Vert
\Phi^{2}_{-\rho}(\tau,t_{0},x)v\right\Vert d\tau \leq
\widetilde{D}(t_{0})\left\Vert
\Phi^{2}_{-\rho}(t,t_{0},x)v\right\Vert,$$ for all $t\geq
t_{0}\geq 0$ and all $(x,v)\in Y$, where $\Phi^{2}_{-\rho}$ is
defined as in Example \ref{ex_shift}. Let functions $M$ and
$\omega$ be given by Definition \ref{def_nedc}. Let us denote
$$c=\int_{0}^{1}e^{-\omega(\tau)}d\tau=\int_{s}^{s+1}e^{-\omega(u-s)}du.$$
Further, for $t\geq s+1$ and $s\geq t_{0}\geq 0$, we obtain
$$c\left\Vert \Phi^{2}_{-\rho}(s,t_{0},x)v\right\Vert=\int_{s}^{s+1}e^{-\omega(u-s)}
\left\Vert\Phi^{2}_{-\rho}(s, t_{0},x)v\right\Vert du\leq$$ $$\leq
\int_{s}^{s+1}M(t_{0})e^{-\omega(u-s)}e^{\omega(u-s)}\left\Vert\Phi^{2}_{-\rho}(u,
t_{0},x)v\right\Vert du\leq$$ $$\leq
M(t_{0})\int_{t_{0}}^{t}\left\Vert\Phi^{2}_{-\rho}(u,
t_{0},x)v\right\Vert du\leq
M(t_{0})\widetilde{D}(t_{0})\left\Vert\Phi^{2}_{-\rho}(t,
t_{0},x)v\right\Vert.$$ We obtain $$\left\Vert
\Phi^{2}_{-\rho}(s,t_{0},x)v\right\Vert\leq
\frac{M(t_{0})\widetilde{D}(t_{0})}{c}\left\Vert\Phi^{2}_{-\rho}(t,
t_{0},x)v\right\Vert,$$ for all $t\geq s\geq t_{0}\geq 0$ with
$t\geq s+1$ and all $(x,v)\in Y.$ Now, for $t\in [s,s+1)$ and
$s\geq t_{0}\geq 0$, we have $$\left\Vert\Phi^{2}_{-\rho}(s,
t_{0},x)v\right\Vert\leq
M(t_{0})e^{\omega(1)}\left\Vert\Phi^{2}_{-\rho}(t,
t_{0},x)v\right\Vert,$$ for all $(x,v)\in Y.$ Finally, we obtain
$$\left\Vert\Phi^{2}_{-\rho}(s,
t_{0},x)v\right\Vert\leq D_{2}(t_{0})\left\Vert\Phi^{2}_{-\rho}(t,
t_{0},x)v\right\Vert,$$ for all $t\geq s\geq t_{0}\geq 0$ and all
$(x,v)\in Y$, where we have denoted
$$D_{2}(u)=M(u)\left[\frac{\widetilde{D}(u)}{c}+e^{\omega(1)}\right], \ u\geq 0.$$
Thus, it follows that $$e^{-(s-t_{0})\rho}\left\Vert\Phi_{2}(s,
t_{0},x)v\right\Vert\leq
D_{2}(t_{0})e^{-(t-t_{0})\rho}\left\Vert\Phi_{2}(t,
t_{0},x)v\right\Vert,$$ which implies
$$\left\Vert\Phi_{2}(s, t_{0},x)v\right\Vert\leq
D_{2}(t_{0})e^{-(t-s)\rho}\left\Vert\Phi_{2}(t,
t_{0},x)v\right\Vert,$$ for all $t\geq s\geq t_{0}\geq 0$ and all
$(x,v)\in Y$, or $$\left\Vert P_{2}(x)v\right\Vert\leq
D_{2}(s)e^{-(t-s)\rho}\left\Vert\Phi_{2}(t,s,x)v\right\Vert,$$ for
all $t\geq s \geq 0$ and all $(x,v)\in Y$.

Hence, the skew-evolution semiflow is exponentially dichotomic,
which ends the proof.
\end{proof}

\vspace{3mm}

\textbf{Acknowledgments.} This work is financially supported from
the Exploratory Research Grant CNCSIS PN II ID 1080 No. 508/2009
of the Romanian Ministry of Education, Research and Innovation.

\vspace{5mm}

\footnotesize{

\vspace{5mm}

\noindent\begin{tabular}[t]{ll}

\textsc{Department of Mathematics and Computer Science}, \\
\textsc{"Aurel Vlaicu"} \textsc{University
of Arad}, \textsc{Romania}\\
\textit{E-mail address:} \texttt{cstoicad@yahoo.com}

\end{tabular}
}

\end{document}